\pdfoutput=1
\RequirePackage{ifpdf}
\ifpdf
\documentclass[pdftex]{sigma}
\else
\documentclass{sigma}
\fi

\DeclareMathOperator{\ev}{ev}
\DeclareMathOperator{\gr}{gr}
\DeclareMathOperator{\spn}{span}

\numberwithin{equation}{section}

\begin{document}

\newcommand{\arXivNumber}{1409.0274}

\allowdisplaybreaks

\renewcommand{\thefootnote}{$\star$}

\renewcommand{\PaperNumber}{110}

\FirstPageHeading

\ShortArticleName{Demazure Modules, Chari--Venkatesh Modules and Fusion Products}

\ArticleName{Demazure Modules, Chari--Venkatesh Modules\\
and Fusion Products\footnote{This paper is a~contribution to the Special Issue on New Directions in Lie Theory.
The full collection is available at
\href{http://www.emis.de/journals/SIGMA/LieTheory2014.html}{http://www.emis.de/journals/SIGMA/LieTheory2014.html}}}

\Author{Bhimarthi RAVINDER}

\AuthorNameForHeading{B.~Ravinder}

\Address{The Institute of Mathematical Sciences, CIT campus, Taramani, Chennai 600113, India}
\Email{\href{mailto:bravinder@imsc.res.in}{bravinder@imsc.res.in}}

\ArticleDates{Received September 11, 2014, in f\/inal form December 01, 2014; Published online December 12, 2014}

\Abstract{Let $\mathfrak{g}$ be a~f\/inite-dimensional complex simple Lie algebra with highest root~$\theta$.
Given two non-negative integers~$m$,~$n$, we prove that the fusion product of~$m$ copies of the level one Demazure
module $D(1,\theta)$ with~$n$ copies of the adjoint representation $\ev_0 V(\theta)$ is independent of the
parameters and we give explicit def\/ining relations.
As a~consequence, for $\mathfrak{g}$ simply laced, we show that the fusion product of a~special family of
Chari--Venkatesh modules is again a~Chari--Venkatesh module.
We also get a~description of the truncated Weyl module associated to a~multiple of~$\theta$.}


\Classification{17B67; 17B10}

\renewcommand{\thefootnote}{\arabic{footnote}} 
\setcounter{footnote}{0}

\section{Introduction}

Let $\mathfrak{g}$ be a~f\/inite-dimensional complex simple Lie algebra with highest root $\theta$. The current algebra
$\mathfrak{g}[t]$ associated to $\mathfrak{g}$ is equal to $\mathfrak{g}\otimes \mathbb{C}[t]$, where $\mathbb{C}[t]$ is
the polynomial ring in one variable.
The degree grading on $\mathbb{C}[t]$ gives a~natural $\mathbb{Z}_{\geq 0}$-grading on $\mathfrak{g}[t]$ and the Lie
bracket is given in the obvious way such that the zeroth grade piece $\mathfrak{g}\otimes 1$ is isomorphic to
$\mathfrak{g}$.
Let $\widehat{\mathfrak{g}}$ be the untwisted af\/f\/ine Lie algebra corresponding to~$\mathfrak{g}$.
In this paper, we shall be concerned with the~$\mathfrak{g}[t]$-stable Demazure modules of integrable highest weight
representations of~$\widehat{\mathfrak{g}}$. The Demazure modules are actually modules for a~Borel subalgebra
$\widehat{\mathfrak{b}}$ of $\widehat{\mathfrak{g}}$. The $\mathfrak{g}[t]$-stable Demazure modules are known to be
indexed by pairs $(l,\lambda)$, where~$l$ is a~positive integer and~$\lambda$ is a~dominant integral weight of~$\mathfrak{g}$ (see~\cite{FoL,Naoi}).
We denote the corresponding module by $D(l, \lambda)$ and call it the level~$l$ Demazure module with highest
weight~$\lambda$; it is in fact a~f\/inite-dimensional graded $\mathfrak{g}[t]$-module.

The study of the category of f\/inite-dimensional graded $\mathfrak{g}[t]$-modules has been of interest in recent years
for variety of reasons.
An important construction in this category is that of the fusion product.
The fusion product of f\/inite-dimensional graded $\mathfrak{g}[t]$-modules~\cite{FL} is by def\/inition, dependent on the
given parameters.
Many people have been working in recent years, to prove the independence of the parameters for the fusion product of
certain~$\mathfrak{g}[t]$-modules, see for instance~\cite{CSVW,CV,FoL,Naoi,V}.
These works mostly considered the fusion product of Demazure modules of the same level and gave explicit def\/ining
relations for them.
We ask the most natural question: Can one give similar results for the fusion product of dif\/ferent level Demazure
modules? In this paper, we answer this question for some important cases; namely we prove (Corollary~\ref{c2}) that the
fusion product of~$m$ copies of the level one Demazure module $D(1, \theta)$ with~$n$ copies of the adjoint
representation $\ev_0 V(\theta)$ is independent of the parameters, and we give explicit def\/ining relations.
We note that $\ev_0 V(\theta)$ may be thought of as a~Demazure module $D(l,\theta)$ of level $l\geq 2$.

More generally, the following is the statement of our main theorem (see Section~\ref{section3} for notation).
\begin{theorem}
\label{MT}
Let $k\geq 1$.
For $0\leq i \leq k$, we have the following:
\begin{enumerate}\itemsep=0pt
\item[$1)$] a~short exact sequence of $\mathfrak{g}[t]$-modules,
\begin{gather*}
0\rightarrow \tau_{2k+1-i} \big(D(1, k\theta)/\big\langle \big(x^{-}_\theta \otimes t^{2k-i}\big) \overline{w}_{k\theta}\big\rangle\big)
\\
\phantom{0}
\xrightarrow{\phi^{-}} D\big(1,(k+1)\theta\big)/\big\langle \big(x^{-}_\theta \otimes t^{2k+2-i}\big) \overline{w}_{(k+1)\theta}\big\rangle
\\
\phantom{0}
\xrightarrow{\phi^+} D\big(1,(k+1)\theta\big)/\big\langle \big(x^{-}_\theta \otimes t^{2k+1-i}\big) \overline{w}_{(k+1)\theta}\big\rangle \rightarrow 0;
\end{gather*}
\item[$2)$] an isomorphism of $\mathfrak{g}[t]$-modules,
\begin{gather*}
D\big(1,(k+1)\theta\big)/\langle \big(x^{-}_\theta \otimes t^{2k+2-i}\big) \overline{w}_{(k+1)\theta}\rangle \cong D(1,\theta)^{* (k+1-i)}* \textup{ev}_0 V(\theta)^{*i}.
\end{gather*}
\end{enumerate}
\end{theorem}

We obtain the following two important corollaries:
\begin{corollary}
\label{c1}
Given $k\geq 1$ and $0\leq i \leq k $, we have the following short exact sequence of $\mathfrak{g}[t]$-modules,
\begin{gather*}
0 \rightarrow \tau_{2k+1-i}\big(D(1,\theta)^{*(k-i)} * \textup{ev}_0 V(\theta)^{*i}\big) \rightarrow D(1,\theta)^{*(k+1-i)}*
\textup{ev}_0 V(\theta)^{*i}
\\
\phantom{0}
\rightarrow D(1,\theta)^{*(k-i)} * \textup{ev}_0 V(\theta)^{*(i+1)} \rightarrow 0.
\end{gather*}
\end{corollary}
\begin{corollary}
\label{c2}
Given $m,n\geq 0$, we have the following isomorphism of $\mathfrak{g}[t]$-modules,
\begin{gather*}
D(1,\theta)^{*m} * \textup{ev}_0 V(\theta)^{*n} \cong D\big(1,(m+n)\theta\big)/\big\langle \big(x^{-}_\theta \otimes t^{2m+n}\big)
\overline{w}_{(m+n)\theta}\big\rangle.
\end{gather*}
\end{corollary}
The Corollary~\ref{c2} generalizes a~result of Feigin (see~\cite[Corollary~2]{F}), where he only considers the case $m=0$.
Theorem~\ref{MT}, Corollaries~\ref{c1} and~\ref{c2} are proved in Section~\ref{section4}.

In~\cite{CV}, Chari and Venkatesh introduced a~large collection of indecomposable graded $\mathfrak{g}[t]$-modules
(which we call Chari--Venkatesh or CV modules) such that all Demazure mo\-du\-les~$D(l, \lambda)$ belong to this collection.
In the case when $\mathfrak{g}$ is simply laced, Theorem~\ref{MT} enables us to obtain (see Theorem~\ref{T2})
interesting exact sequences between CV modules and to show that the fusion product of a~special family of CV modules is
again a~CV module.
Theorem~\ref{T2} generalizes results of Chari and Venkatesh (see~\cite[\S~6]{CV}), where they only consider the case $\mathfrak{g}=\mathfrak{sl}_2$.

For $n\geq 1$, let $\mathcal{A}_n=\mathbb{C}[t]/(t^n)$ be the truncated algebra.
We consider for $k\geq1$ the local Weyl modules $W_{\mathcal{A}_n}(k\theta)$ for the truncated current algebra
$\mathfrak{g}\otimes\mathcal{A}_n$.
These modules are known to be f\/inite-dimensional, but they are still far from being well understood; even their
dimensions are not known.
As a~consequence of Theorem~\ref{MT}, we are able to obtain the following description of truncated Weyl modules in terms
of local Weyl modules $W(k\theta)$, $k\geq 1$, for the current algebra~$\mathfrak{g}[t]$.
The latter modules $W(k\theta)$ are very well understood.
\begin{corollary}
\label{truncated}
Assume that $\mathfrak{g}$ is simply laced.
Given $k,n\geq 1$, we have the following isomorphism of $\mathfrak{g}[t]$-modules,
\begin{gather*}
W_{\mathcal{A}_n}(k\theta) \cong
\begin{cases}
W(\theta)^{*(n-k)} * \textup{ev}_0 V(\theta)^{*(2k-n)}, & k\leq n < 2k,
\\
W(k\theta), & n\geq 2k.
\end{cases}
\end{gather*}
\end{corollary}

 The Corollary~\ref{truncated} is proved in Section~\ref{section5}.

\section{Preliminaries}\label{section2}

Throughout the paper, $\mathbb{C}$ denote the f\/ield of complex numbers, $\mathbb{Z}$ the set of integers,
$\mathbb{Z}_{\geq 0}$ the set of non-negative integers, $\mathbb{N}$ the set of positive integers and $\mathbb{C}[t]$
the polynomial ring in an indeterminate~$t$.

{\bf 2.1.}~Let $\mathfrak{a}$ be a~complex Lie algebra, $\mathbf{U}(\mathfrak{a})$ the corresponding universal
enveloping algebra.
The current algebra associated to $\mathfrak{a}$ is denoted by $\mathfrak{a}[t]$ and def\/ined as $\mathfrak{a} \otimes
\mathbb{C}[t]$, with the Lie bracket
\begin{gather*}
[a \otimes t^r, b \otimes t^s]=[a, b]\otimes t^{r+s},
\qquad
\text{for all}
\quad
a, b \in \mathfrak{a}
\quad
\text{and}
\quad
r, s \in \mathbb{Z}_{\geq 0}.
\end{gather*}
We let $\mathfrak{a}[t]_{+}$ be the ideal $\mathfrak{a}\otimes t\mathbb{C}[t]$. The degree grading on $\mathbb{C}[t]$
gives a~natural $\mathbb{Z}_{\geq 0}$-grading on~$\mathbf{U}(\mathfrak{a}[t])$ and the subspace of grade~$s$ is given~by
\begin{gather*}
\mathbf{U}(\mathfrak{a}[t])[s]= \spn \Big\{(a_1\otimes t^{r_1})\cdots (a_k\otimes t^{r_k}):k\geq 1,\,
a_i\in\mathfrak{a},\, r_i\in \mathbb{Z}_{\geq 0}, \sum r_i=s\Big\},
\quad
\forall\, s\in \mathbb{N},
\end{gather*}
and the subspace of grade zero $\mathbf{U}(\mathfrak{a}[t])[0]=\mathbf{U}(\mathfrak{a})$.

{\bf 2.2.}~Let $\mathfrak{g}$ be a~f\/inite-dimensional complex simple Lie algebra, with Cartan subalgebra~$\mathfrak{h}$.
Let~$R$ (resp.~$R^+$) be the set of roots (resp.\
positive roots) of $\mathfrak{g}$ with respect to $\mathfrak{h}$ and $\theta \in R^+$ be the highest root in~$R$.
There is a~non-degenerate, symmetric, Weyl group invariant bilinear form $(\cdot|\cdot)$ on $\mathfrak{h}^*$, which we assume to
be normalized so that the square length of a~long root is two.
For $\alpha\in R$, $\alpha^{\vee}\in\mathfrak{h}$ denotes the corresponding co-root and we set
$d_{\alpha}=2/(\alpha|\alpha)$. For $\alpha\in R$, let $\mathfrak{g}_{\alpha}$ be the corresponding root space of
$\mathfrak{g}$ and f\/ix non-zero elements $x^{\pm}_{\alpha}\in \mathfrak{g}_{\pm\alpha}$ such that $[x^{+}_{\alpha},x^{-}_{\alpha}]=\alpha^{\vee}$.
We set $\mathfrak{n}^{\pm}=\oplus_{\alpha\in R^{+}}  \mathfrak{g}_{\pm\alpha}$.

Let $P^{+}$ be the set of dominant integral weights of $\mathfrak{g}$.
For $\lambda\in P^+$, $V(\lambda)$ be the corresponding f\/inite-dimensional irreducible $\mathfrak{g}$-module generated~by
an element $v_{\lambda}$ with the following def\/ining relations:
\begin{gather*}
x^{+}_{\alpha}v_\lambda=0,
\qquad
h v_\lambda = \langle \lambda, h \rangle v_\lambda,
\qquad
(x^{-}_{\alpha})^{\langle \lambda, \alpha^{\vee} \rangle +1}  v_\lambda=0,
\qquad
\text{for all}
\quad
\alpha\in R^+,
\quad
h\in\mathfrak{h}.
\end{gather*}

{\bf 2.3.}~A graded $\mathfrak{g}[t]$-module is a~$\mathbb{Z}$-graded vector space
\begin{gather*}
V=\bigoplus_{r\in\mathbb{Z}} V[r] \qquad
\text{such that}
\quad
(x\otimes t^s) V[r]\subset V[r+s],
\quad
x\in\mathfrak{g},
\quad
r\in\mathbb{Z},
\quad
s\in\mathbb{Z}_{\geq0}.
\end{gather*}
For $\mu\in\mathfrak{h}^*$, an element~$v$ of a~graded $\mathfrak{g}[t]$-module~$V$ is said to be of weight~$\mu$, if
$(h\otimes 1)v=\langle \mu, h\rangle v$ for all $h\in \mathfrak{h}$.
We def\/ine a~morphism between two graded $\mathfrak{g}[t]$-modules as a~degree zero morphism of
$\mathfrak{g}[t]$-modules.
For $r\in\mathbb{Z}$, let $\tau_r$ be the grade shift operator: if~$V$ is a~graded $\mathfrak{g}[t]$-module then~$\tau_r
V$ is the graded $\mathfrak{g}[t]$-module with the graded pieces shifted uniformly by~$r$ and the action of~$\mathfrak{g}[t]$ remains unchanged.
For any graded $\mathfrak{g}[t]$-module~$V$ and a~subset~$S$ of~$V$, $\langle S\rangle$ denotes the submodule of~$V$ generated by~$S$.
For $\lambda\in P^+$, $\ev_0 V(\lambda)$ be the irreducible graded $\mathfrak{g}[t]$-module such
that $\ev_0 V(\lambda)[0]\cong_{\mathfrak{g}} V(\lambda)$ and $\ev_0 V(\lambda)[r]=0$
$\forall\, r\in \mathbb{N}$.
In particular, $\mathfrak{g}[t]_{+}(\ev_0 V(\lambda))=0$.

{\bf 2.4.}~For $r,s\in\mathbb{Z}_{\geq0}$, we denote
\begin{gather*}
\mathbf{S}(r,s)=\bigg\{(b_p)_{p\geq0}: b_p\in\mathbb{Z}_{\geq0},
\;
\sum\limits_{p\geq 0}b_p =r,
\;
\sum\limits_{p\geq0}pb_p=s\bigg\}.
\end{gather*}
For $\alpha\in R^{+}$ and $r,s\in\mathbb{Z}_{\geq0}$, we def\/ine an element
$\mathbf{x}^{-}_{\alpha}(r,s)\in\mathbf{U}(\mathfrak{g}[t])[s]$~by
\begin{gather*}
\mathbf{x}^{-}_{\alpha}(r,s)=\sum\limits_{(b_p)\in\mathbf{S}(r,s)} (x^{-}_{\alpha} \otimes 1)^{(b_0)}(x^{-}_{\alpha}
\otimes t)^{(b_1)}\cdots(x^{-}_{\alpha} \otimes t^s)^{(b_s)},
\end{gather*}
where for any non-negative integer~$b$ and any $x\in\mathfrak{g}[t]$, we understand
$x^{(b)}=x^b/b!$.

The following was proved in~\cite{G} (see also~\cite[Lemma 2.3]{CV}).
\begin{lemma}\label{garland}
Given $s\in\mathbb{N}$, $r\in\mathbb{Z}_{\geq0}$ and $\alpha\in R^{+}$, we have
\begin{gather*}
(x^{+}_{\alpha} \otimes t)^{(s)}(x^{-}_{\alpha} \otimes
1)^{(s+r)}-(-1)^s\mathbf{x}^{-}_{\alpha}(r,s)\in\mathbf{U}(\mathfrak{g}[t])\mathfrak{n}^{+}[t]\bigoplus
\mathbf{U}(\mathfrak{n}^{-}[t])\mathfrak{h}[t]_{+}.
\end{gather*}
\end{lemma}

\section{Weyl, Demazure modules and fusion product}\label{section3}

In this section, we recall the def\/initions of local Weyl modules, level one Demazure modules and fusion products.

\subsection{Weyl module}\label{section3.1}

The def\/inition of the local Weyl module was given originally in~\cite{CP}, later in~\cite{CFK} and~\cite{FL}.
\begin{definition}
Given $\lambda \in P^+$, the local Weyl module $W(\lambda)$ is the cyclic $\mathfrak{g}[t]$-module generated by an
element $w_\lambda$, with following def\/ining relations:
\begin{gather}
\mathfrak{n}^{+}[t]  w_\lambda=0,
\qquad
(h \otimes t^s)  w_\lambda = \langle \lambda, h \rangle \delta_{s,0} w_\lambda=0,
\qquad
s\geq0,
\qquad
h\in\mathfrak{h},
\nonumber
\\
(x^{-}_\alpha \otimes 1)^{\langle \lambda,
\alpha^{\vee} \rangle + 1}w_\lambda=0,
\qquad
\alpha \in R^{+}.
\label{w2}
\end{gather}
\end{definition}
\noindent We note that the relation~\eqref{w2} implies
\begin{gather}
\label{w2'}
\big(x^{-}_\alpha \otimes t^{\langle \lambda,
\alpha^{\vee} \rangle}\big)  w_\lambda=0,
\qquad
\alpha \in R^{+},
\end{gather}
which is easy to see from Lemma~\ref{garland}.
We set the grade of $w_\lambda$ to be zero; then $W(\lambda)$ becomes a~$\mathbb{Z}_{\geq 0}$-graded module with
\begin{gather*}
W(\lambda)[0] \cong_{\mathfrak{g}} V(\lambda).
\end{gather*}
Moreover, $\textup{ev}_0 V(\lambda)$ is the unique graded irreducible quotient of $W(\lambda)$.

We now specialize to the case $\lambda\in \mathbb{N}\theta$, and obtain some further useful relations that hold in~$W(\lambda)$.
\begin{lemma}
\label{g}
Let $k\in\mathbb{N}$.
The following relations hold in the local Weyl module $W((k+1)\theta)$:
\begin{enumerate}\itemsep=0pt
\item[$1)$] $(x^{-}_{\theta}\otimes 1)^{2k+1}\big(x^{-}_{\theta}\otimes t^{2k+1-i}\big)w_{(k+1)\theta}=0,
\qquad
\forall\,
0\leq i \leq k$;
\item[$2)$] $(x^{-}_{\theta}\otimes t^m)(x^{-}_{\theta} \otimes t^{m+1})  w_{(k+1)\theta} \in  \big\langle \big(x^{-}_{\theta}
\otimes t^{m+2}\big) w_{(k+1)\theta} \rangle,
\qquad
\forall\,
m\geq k$.
\end{enumerate}
\end{lemma}

\begin{proof}
To prove part (1), consider $(x^{+}_{\theta}\!\otimes t^{2k{+}1{-}i}) (x^{-}_{\theta}\otimes 1)^{2k{+}3}w_{(k{+}1)\theta}$.
Since $(x^{+}_{\theta}\!\otimes t^{2k{+}1{-}i})w_{(k{+}1)\theta}$ $=0$, we get
\begin{gather*}
\big(x^{+}_{\theta}\otimes t^{2k+1-i}\big) (x^{-}_{\theta}\otimes 1)^{2k+3}w_{(k+1)\theta}
=\big[x^{+}_{\theta}\otimes t^{2k+1-i}, (x^{-}_{\theta}\otimes 1)^{2k+3}\big] w_{(k+1)\theta}
\\
\qquad{}
=\sum\limits_{j=1}^{2k+3}  (x^{-}_{\theta}\otimes 1)^{j-1}\big(\theta^{\vee}\otimes t^{2k+1-i}\big)(x^{-}_{\theta}\otimes
1)^{2k+3-j}w_{(k+1)\theta}.
\end{gather*}
Since $(\theta^{\vee}\otimes t^{2k+1-i})w_{(k+1)\theta}=0$, we may replace $(\theta^{\vee}\otimes
t^{2k+1-i})(x^{-}_{\theta}\otimes 1)^{2k+3-j}$~by
\begin{gather*}
\big[\theta^{\vee}\otimes t^{2k+1-i}, (x^{-}_{\theta}\otimes 1)^{2k+3-j}\big]
=(-2)(2k+3-j)(x^{-}_{\theta}\otimes 1)^{2k+2-j}\big(x^{-}_{\theta}\otimes t^{2k+1-i}\big).
\end{gather*}
After simplifying, we get
\begin{gather*}
\big(x^{+}_{\theta}\otimes t^{2k+1-i}\big) (x^{-}_{\theta}\otimes 1)^{2k+3}w_{(k+1)\theta}
\\
\qquad{}
=(-1)(2k+2)(2k+3)(x^{-}_{\theta}\otimes 1)^{2k+1}\big(x^{-}_{\theta}\otimes t^{2k+1-i}\big)w_{(k+1)\theta}.
\end{gather*}
Now, using $(x^{-}_{\theta}\otimes 1)^{2k+3}w_{(k+1)\theta}=0$ in $W((k+1)\theta)$, completes the proof of part~(1).
Part (2) follows easily by putting $r=2$, $s=2m+1$ and $\alpha=\theta$ in Lemma~\ref{garland}, and using the fact that
$(x^{-}_{\theta}\otimes 1)^{2m+3}w_{(k+1)\theta}=0$, $\forall\,
m\geq k$ by~\eqref{w2}.
\end{proof}

\subsection{Level one Demazure module}\label{section3.2}

Let $\lambda\in P^{+}$ and $\alpha\in R^{+}$ with $\langle \lambda, \alpha^{\vee}\rangle > 0$.
Let $s_\alpha, m_\alpha \in \mathbb{N}$ be the unique positive integers such that
\begin{gather*}
\langle \lambda, \alpha^{\vee}\rangle = (s_\alpha-1)d_\alpha + m_\alpha,
\qquad
0<m_\alpha\leq d_\alpha.
\end{gather*}
If $\langle \lambda, \alpha^{\vee}\rangle = 0$, set $s_\alpha=0=m_\alpha$.
We take the following as a~def\/inition of the level one Demazure module.
\begin{definition}
\textup{(see~\cite[Corollary 3.5]{CV})} The level one Demazure module $D(1,\lambda)$ is the graded quotient of
$W(\lambda)$ by the submodule generated by the union of the following two sets:
\begin{gather}
\big\{(x^{-}_{\alpha} \otimes t^{s_{\alpha}})  w_{\lambda}: \alpha\in R^{+}~\textup{such that}~d_{\alpha} > 1\big\},
\label{dm1}
\\
\big\{\big(x^{-}_{\alpha} \otimes t^{s_{\alpha}-1}\big)^2  w_{\lambda}: \alpha\in R^{+}~\textup{such that}~d_\alpha =3~\textup{and}~m_\alpha=1\big\}.
\label{dm2}
\end{gather}
In particular, for $\mathfrak{g}$ simply laced, $D(1,\lambda)\cong_{\mathfrak{g}[t]} W(\lambda)$. We denote~by
$\overline{w}_\lambda$, the image of $w_\lambda$ in~$D(1,\lambda)$.
\end{definition}

The following proposition gives explicit def\/ining relations for $D(1,k\theta)$.
\begin{proposition}
\label{WvsD}
Given $k\geq 1$, the level~$1$ Demazure module $D(1,k\theta)$ is the graded $\mathfrak{g}[t]$-module generated by an
element $\overline{w}_{k\theta}$, with the following defining relations:
\begin{gather*}
\mathfrak{n}^{+}[t]\, \overline{w}_{k\theta}=0,
\qquad
(h \otimes t^s) \overline{w}_{k\theta} = \langle k\theta, h\rangle \delta_{s,0} \overline{w}_{k\theta},
\qquad
s\geq0,
\quad
h\in\mathfrak{h},
\\
(x^{-}_{\alpha}\otimes 1) \overline{w}_{k\theta}=0,
\qquad
\alpha\in R^+,
\quad
(\theta|\alpha)=0,
\\
(x^{-}_{\alpha} \otimes 1)^{kd_{\alpha}+1} \overline{w}_{k\theta}=0,
\qquad
\big(x^{-}_{\alpha} \otimes t^k\big) \overline{w}_{k\theta}=0,
\qquad
\alpha\in R^+,
\quad
(\theta|\alpha)=1,
\\
(x^{-}_{\theta} \otimes 1)^{2k+1} \overline{w}_{k\theta}=0.
\end{gather*}
\end{proposition}
\begin{proof}
Observe that, from the abstract theory of root systems $(\theta|\alpha)= 0$ or~1,
$\forall\, \alpha \in R^{+}\setminus \{\theta\}$.
This implies that $\langle k\theta, \alpha^{\vee}\rangle = 0$ or $kd_{\alpha}$, $\forall\, \alpha \in R^{+}\setminus \{\theta\}$.
Hence the relations~\eqref{dm2} do not occur in $D(1,k\theta)$ and the relations~\eqref{dm1} are
\begin{gather*}
\big(x^{-}_{\alpha} \otimes t^k\big) \overline{w}_{k\theta}=0,
\qquad
\alpha\in R^+,
\quad
\alpha~~\text{short},
\quad
(\theta|\alpha)=1.
\end{gather*}
For a~long root $\alpha\in R^+$ with $(\theta|\alpha)=1$, by~\eqref{w2'} it follows that $(x^{-}_{\alpha} \otimes t^k) \overline{w}_{k\theta}=0$.
Now the other relations are precisely the def\/ining relations of $W(k\theta)$. This proves Proposition~\ref{WvsD}.
\end{proof}
We record below a~well-known fact, for later use:
\begin{gather*}
D(1, \theta)\cong_{\mathfrak{g}} V(\theta) \oplus \mathbb{C}.
\end{gather*}
In particular,
\begin{gather}
\label{dimd}
\dim D(1, \theta)= \dim V(\theta)+1.
\end{gather}

The following is a~crucial lemma, which we use in proving Theorem~\ref{MT}.
\begin{lemma}
\label{crucial}
Let $k\geq 1$ and $0\leq i \leq k$. The following relations hold in the module $D(1,(k+1)\theta)$:
\begin{enumerate}\itemsep=0pt
\item[$1)$] $(x^{-}_{\alpha} \otimes 1)^{kd_{\alpha}+1}\big(x^{-}_{\theta} \otimes t^{2k+1-i}\big) \overline{w}_{(k+1)\theta}=0$,
$\forall\, \alpha \in R^{+}$, $(\theta|\alpha)=1$;
\item[$2)$] $\big(x^{-}_{\alpha} \otimes t^{k}\big) (x^{-}_{\theta} \otimes t^{2k+1-i}) \overline{w}_{(k+1)\theta} \in
\big\langle \big(x^{-}_{\theta} \otimes t^{2k+2-i}\big) \overline{w}_{(k+1)\theta} \big\rangle$,
$\forall\, \alpha \in R^{+}$, $(\theta|\alpha)=1$;
\item[$3)$] $\big(x^{-}_{\theta} \otimes t^{2k-i}\big) \big(x^{-}_{\theta} \otimes t^{2k+1-i}\big) \overline{w}_{(k+1)\theta}
\in
\big\langle \big(x^{-}_{\theta} \otimes t^{2k+2-i}\big) \overline{w}_{(k+1)\theta} \big\rangle$.
\end{enumerate}
\end{lemma}
\begin{proof}
Let $ \alpha \in R^{+} \textup{with} (\theta|\alpha)=1$. This implies that $\theta-\alpha$ is also a~root of
$\mathfrak{g}$ and $(\theta|\theta-\alpha)=1$. We now prove part (1).
Observe that, $(x^{-}_{\theta} \otimes t^{2k+1-i}) \overline{w}_{(k+1)\theta}$ is an element of weight $k\theta$.
Further $(x^{+}_{\alpha} \otimes 1) (x^{-}_{\theta} \otimes t^{2k+1-i}) \overline{w}_{(k+1)\theta}=0$, since
$(x^{+}_{\alpha} \otimes 1)\overline{w}_{(k+1)\theta}=0$ and $(x^{-}_{\theta-\alpha} \otimes t^{2k+1-i}) \overline{w}_{(k+1)\theta}=0$,
for all $0\leq i \leq k$.
Considering the copy of $\mathfrak{sl}_2$ spanned by $x^{+}_{\alpha} \otimes 1$, $x^{-}_{\alpha} \otimes 1$,
$\alpha^{\vee}\otimes 1$, we obtain part (1) by standard $\mathfrak{sl}_2$ arguments.
We now prove part (2).
Putting $r=2$, $s=(3k+1-i)$ and $\alpha=\theta$ in Lemma~\ref{garland}, we get
\begin{gather}
\big(x^{-}_{\theta} \otimes t^{k}\big) \big(x^{-}_{\theta} \otimes t^{2k+1-i}\big) \overline{w}_{(k+1)\theta} +
\sum\limits_{\substack{k+1\leq p\leq q \leq 2k-i\\p+q= 3k+1-i}}
\frac{1}{(2\delta_{p, q})!} \big(x^{-}_{\theta} \otimes t^{p}\big) \big(x^{-}_{\theta} \otimes t^{q}\big)
\overline{w}_{(k+1)\theta}
\nonumber\\
\qquad{}
\in
\big\langle\big(x^{-}_{\theta}\otimes t^{2k+2-i}\big) \overline{w}_{(k+1)\theta} \big\rangle, \label{e}
\end{gather}
since $(x^{-}_{\theta}\otimes 1)^{3k+3-i}\overline{w}_{(k+1)\theta}=0$, $\forall\, 0\leq i \leq k$.
Now we act on both sides of~\eqref{e} by $x^{+}_{\theta-\alpha}$ and use the relation $(x^{-}_{\alpha} \otimes t^r) \overline{w}_{(k+1)\theta}=0$,
for all $r\geq(k+1)$, which gives part~(2).
Part~(3) is immediate from the part~(2) of Lemma~\ref{g}.
\end{proof}

\subsection{Fusion product}\label{section3.3}

In this subsection, we recall the def\/inition of the fusion product of f\/inite-dimensional graded cyclic
$\mathfrak{g}[t]$-modules given in~\cite{FL} and give some elementary properties.

For a~cyclic $\mathfrak{g}[t]$-module~$V$ generated by~$v$, we def\/ine a~f\/iltration $F^{r}V$, $r\in\mathbb{Z}_{\geq 0}$~by
\begin{gather*}
F^{r}V=\sum\limits_{0\leq s \leq r} \mathbf{U}(\mathfrak{g}[t])[s] v.
\end{gather*}
We say $F^{-1}V$ is the zero space.
The associated graded space $\gr V=\bigoplus_{r\geq 0} F^{r}V/ F^{r-1}V $ naturally becomes a~cyclic
$\mathfrak{g}[t]$-module generated by $v+F^{-1}V$, with action given by
\begin{gather*}
(x\otimes t^s)\big(w+F^{r-1}V\big):= (x\otimes t^s)w+F^{r+s-1}V,
\qquad
\forall\,
x\in \mathfrak{g},
\quad
w\in F^{r}V, \quad
r, s\in \mathbb{Z}_{\geq0}.
\end{gather*}
Observe that, $\gr V \cong V$ as $\mathfrak{g}$-modules.

The following lemma will be useful.
\begin{lemma}
\label{f1}
Let~$V$ be a~cyclic $\mathfrak{g}[t]$-module.
For $r,s\in \mathbb{Z}_{\geq0}$, the following equality holds in the quotient space $F^{r+s}V/ F^{r+s-1}V$.
\begin{gather*}
(x\otimes t^s)\big(w+F^{r-1}V\big)=\big((x\otimes (t-a_1)\cdots (t-a_s))w\big) + F^{r+s-1}V,
\end{gather*}
for all $a_1,\dots,a_s\in\mathbb{C}$,
$x\in \mathfrak{g}$, $w\in F^{r}V$.
\end{lemma}

Given a~$\mathfrak{g}[t]$-module~$V$ and $z\in\mathbb{C}$, we def\/ine an another $\mathfrak{g}[t]$-module action on~$V$
as follows:
\begin{gather*}
(x\otimes t^s)v=\big(x\otimes (t+z)^s\big)v,
\qquad
x\in\mathfrak{g},
\qquad
v\in V,
\qquad
s\in \mathbb{Z}_{\geq0}.
\end{gather*}
We denote this new module by $V^z$.

Let $V_i$ be a~f\/inite-dimensional cyclic graded $\mathfrak{g}[t]$-module generated by $v_i$, for $1\leq i\leq m$, and
let $z_1,\dots,z_m$ be distinct complex numbers.
We denote~by
\begin{gather*}
\mathbf{V}={V_1}^{z_1}\otimes\dots\otimes {V_m}^{z_m},
\end{gather*}
the corresponding tensor product of $\mathfrak{g}[t]$-modules.
It is easily checked (see~\cite[Proposition 1.4]{FL}) that $\mathbf{V}$ is a~cyclic $\mathfrak{g}[t]$-module generated
by $v_1\otimes\dots\otimes v_m$.
The associated graded space $\gr \mathbf{V}$ is called the fusion product of $V_1,\dots,V_m$ w.r.t.\ the parameters
$z_1,\dots,z_m$, and is denoted by ${V_1}^{z_1}*\cdots*{V_m}^{z_m}$.
We denote $v_1*\cdots*v_m=(v_1\otimes\cdots\otimes v_m)+F^{-1}\mathbf{V}$, a~generator of $\gr \mathbf{V}$.
For ease of notation we mostly, just write $V_1*\cdots*V_m$ for ${V_1}^{z_1}*\cdots*{V_m}^{z_m}$.
But unless explicitly stated, it is assumed that the fusion product does depend on these parameters.

The following lemma will be needed later.
\begin{lemma}
\label{f2}
Given $1\leq i\leq m$, let $V_i$ be a~finite-dimensional cyclic graded $\mathfrak{g}[t]$-module generated by $v_i$, and
$s_i\in\mathbb{Z}_{\geq0}$.
Let $x\in\mathfrak{g}$.
If $(x\otimes t^{s_i})v_i=0$, $\forall\, 1\leq i\leq m$ then $(x\otimes t^{s_1+\dots+s_m}) v_1*\dots*v_m=0$.
\end{lemma}
\begin{proof}
Let $z_1,\dots,z_m$ be distinct complex numbers and let $\mathbf{V}$ as above.
By using Lemma~\ref{f1}, we get the following equality in $\gr \mathbf{V}$,
\begin{gather*}
\big(x\otimes t^{s_1+\cdots+s_m}\big) \big((v_1\otimes\cdots\otimes v_m)+F^{-1}\mathbf{V}\big)
\\
\qquad{}
=\big(\big(x\otimes (t-z_1)^{s_1}\cdots(t-z_m)^{s_m}\big)v_1\otimes\dots\otimes v_m\big) + F^{s_1+\cdots+s_m-1}\mathbf{V}.
\end{gather*}
Now the proof follows by the def\/inition of the fusion product.
\end{proof}

\section{Proof of the main theorem}\label{section4}

In this section, we prove the existence of maps $\phi^{+}$ and $\phi^{-}$ from Theorem~\ref{MT} and then prove our main
theorem (Theorem~\ref{MT}).

{\bf 4.1.}~Given $k\geq 1$ and $0\leq i \leq k $, we denote~by
\begin{gather*}
\mathbf{V}_{i,k} = D(1, k\theta)/\big\langle \big(x^{-}_\theta \otimes t^{2k-i}\big) \overline{w}_{k\theta} \big\rangle,
\end{gather*}
and let $\overline{v}_{i,k}$ be the image of $\overline{w}_{k\theta}$ in $\mathbf{V}_{i,k}$.

Using Proposition~\ref{WvsD}, $\mathbf{V}_{i,k}$ is the cyclic graded $\mathfrak{g}[t]$-module generated by the element
$\overline{v}_{i,k}$, with the following def\/ining relations:
\begin{gather}
(x^{+}_{\alpha}\otimes t^s)\overline{v}_{i,k}=0,
\qquad
s\geq 0,\quad \alpha\in R^+,
\label{r1}
\\
(h \otimes t^s) \overline{v}_{i,k} = \langle k\theta, h\rangle \delta_{s,0} \overline{v}_{i,k},
\qquad
s\geq0,
\quad
h\in\mathfrak{h},
\label{r2}
\\
(x^{-}_{\alpha}\otimes 1) \overline{v}_{i,k}=0,
\qquad
\alpha\in R^+,
\quad
(\theta|\alpha)=0,
\label{r3}
\\
(x^{-}_{\alpha}\otimes 1)^{kd_{\alpha}+1} \overline{v}_{i,k}=0,
\qquad
\big(x^{-}_{\alpha} \otimes t^k\big) \overline{v}_{i,k}=0,
\qquad
\alpha\in R^+,
\quad
(\theta|\alpha)=1,
\label{r4}
\\
(x^{-}_{\theta}\otimes 1)^{2k+1} \overline{v}_{i,k}=0,
\qquad
\big(x^{-}_{\theta} \otimes t^{2k-i}\big) \overline{v}_{i,k}=0.
\label{r5}
\end{gather}
The existence of $\phi^{+}$ is trivial, which we record below.
\begin{proposition}
The map $\phi^{+} : \mathbf{V}_{i,k+1} \rightarrow \mathbf{V}_{i+1,k+1}$ which takes $\overline{v}_{i,k+1}
\rightarrow \overline{v}_{i+1,k+1}$ is a~surjective morphism of $\mathfrak{g}[t]$-modules with $\ker \phi^{+} = \big\langle
\big({x^{-}_\theta} \otimes t^{2k+1-i}\big) \overline{v}_{i,k+1}\big\rangle$.
\end{proposition}

Now we prove the existence of $\phi^{-}$ in the following proposition.
\begin{proposition}
There exist a~surjective morphism of $\mathfrak{g}[t]$-modules $\phi^{-} : \tau_{2k+1-i} \mathbf{V}_{i,k}
\rightarrow \ker \phi^{+}$,
such that $\phi^{-}(\overline{v}_{i,k}) = \big(x^{-}_\theta \otimes t^{2k+1-i}\big) \overline{v}_{i,k+1}$.
\end{proposition}

\begin{proof}
We only need to show that $\phi^{-}(\overline{v}_{i,k})$ satisf\/ies the def\/ining relations of $\mathbf{V}_{i,k}$. We
start with the relation~\eqref{r1}.
First, for $\alpha=\theta$ it is clear.
Let $\alpha \in R^{+}\setminus\{\theta\}$; if $(\theta|\alpha)=0$ then also it is clear.
If $(\theta|\alpha)=1$ then $(\theta-\alpha)\in R^{+}\setminus\{\theta\}$ and $(\theta|\theta-\alpha)=1$, now it is
clear from the relations $(x^{-}_{\theta-\alpha}\otimes t^r)\overline{v}_{i,k+1}=0$ for all $r\geq (k+1)$ in
$\mathbf{V}_{i,k+1}$.
The relations~\eqref{r2},~\eqref{r3} are trivially satisf\/ied by $\phi^{-}(\overline{v}_{i,k})$.
Finally the last two relations~\eqref{r4},~\eqref{r5} are also satisf\/ied by $\phi^{-}(\overline{v}_{i,k})$; in fact
these are exactly the statements of Lemmas~\ref{g} and~\ref{crucial}.
\end{proof}

{\bf 4.2.}~The existence of the surjective maps $\phi^{+}$ and $\phi^{-}$, give the following:
\begin{gather}
\label{diml}
\dim \mathbf{V}_{i,k+1} \leq \dim \mathbf{V}_{i,k} + \dim \mathbf{V}_{i+1,k+1}.
\end{gather}

The following proposition helps in proving the reverse inequality.
\begin{proposition}
The map $\psi :\mathbf{V}_{i,k+1}\rightarrow D(1,\theta)^{* (k+1-i)} * \textup{ev}_0 V(\theta)^{*i}$ such that
$\psi(\overline{v}_{i,k+1}) = \overline{w}_{\theta}^{* (k+1-i)} * v_{\theta}^{*i}$ is well-defined and surjective
morphism of $\mathfrak{g}[t]$-modules.
In particular,
\begin{gather}
\label{dimg}
\dim \mathbf{V}_{i,k+1}\geq (\dim D(1, \theta))^{k+1-i} (\dim V(\theta))^{i}.
\end{gather}
\end{proposition}
\begin{proof}
We only need to show that $\psi(\overline{v}_{i,k+1})$ satisf\/ies the def\/ining relations of $\mathbf{V}_{i,k+1}$. But
they follow easily from the following relations:
\begin{gather*}
((h \otimes 1)-\langle (k+1)\theta, h\rangle)\big(\overline{w}_{\theta}^{\otimes (k+1-i)} \otimes v_{\theta}^{\otimes i}\big)=0,
\qquad
\forall\, h\in\mathfrak{h},
\\
(x^{-}_{\alpha}\otimes 1)^{\langle (k+1)\theta,
\qquad
\alpha^{\vee} \rangle + 1}\big(\overline{w}_{\theta}^{\otimes (k+1-i)}\otimes v_{\theta}^{\otimes i}\big)=0,
\qquad
\forall\, \alpha\in R^{+}
\end{gather*}
(which holds in $D(1,\theta)^{\otimes (k+1-i)} \otimes \textup{ev}_0 V(\theta)^{\otimes i}$) and $(h \otimes
t^s)\psi(\overline{v}_{i,k+1})=0, \forall\, s\geq 1$, $h\in\mathfrak{h}$
(which holds in $D(1,\theta)^{* (k+1-i)}*\textup{ev}_0 V(\theta)^{*i}$).
Further from Lemma~\ref{f2}, by using the relations
\begin{gather*}
(x^{+}_{\alpha}\otimes t^s)\overline{w}_{\theta}=0=(x^{+}_{\alpha}\otimes t^s)v_{\theta},
\qquad
\forall\, s\geq 0, \quad \alpha\in R^{+},
\\
(x^{-}_{\alpha}\otimes t)\overline{w}_{\theta}=\big(x^{-}_{\theta}\otimes
t^{2}\big)\overline{w}_{\theta}=0=(x^{-}_{\theta}\otimes t)v_{\theta}=(x^{-}_{\alpha}\otimes t)v_{\theta},
\qquad
\forall\, \alpha \in R^{+}\setminus\{\theta\},
\end{gather*}
which holds in $D(1,\theta)$ and $\textup{ev}_0 V(\theta)$.
\end{proof}

We record below a~result from~\cite{F} and use this in proving our main theorem.
\begin{proposition}\textup{\cite[Corollary 2]{F}}
\label{F}
Given $k\geq 1$, the following is an isomorphism of $\mathfrak{g}[t]$-modules,
\begin{gather*}
\textup{ev}_0 V(\theta)^{* k} \cong D(1, k\theta)/\big\langle \big({x^{-}_\theta} \otimes t^{k}\big) \overline{w}_{k\theta}\big\rangle.
\end{gather*}
\end{proposition}

{\bf 4.3.}~We now prove Theorem~\ref{MT}, proceeding by induction on~$k$.
First, for $k=1$, we prove Theorem~\ref{MT} for $0\leq i \leq 1$.
Let $i=1$, observe that $\mathbf{V}_{1,1} \cong_{\mathfrak{g}[t]} \ev_{0} V(\theta)$.
Using Proposition~\ref{F},~\eqref{dimd},~\eqref{diml} and~\eqref{dimg} this case follows.
Let $i=0$, now observe that $\mathbf{V}_{0,1} \cong_{\mathfrak{g}[t]} D(1,\theta)$.
Using part (2) of Theorem~\ref{MT} for $i=1$ and $k=1$,~\eqref{dimd},~\eqref{diml} and~\eqref{dimg} this case also
follows.
Now let $k\geq2$, and assume Theorem~\ref{MT} holds for $(k-1)$.
We prove the assertion for~$k$, proceeding by induction on~$i$.
For $i=k$, it follows from Proposition~\ref{F},~\eqref{dimd},~\eqref{diml} and~\eqref{dimg}.
Now let $i\leq (k-1)$, and assume Theorem~\ref{MT} holds for $(i+1)$.
We now prove for~$i$.
Using part (2) of Theorem~\ref{MT}, for $(i+1)$ and~$k$, also for~$i$ and $(k-1)$, and~\eqref{diml}, we get
\begin{gather*}
\dim \mathbf{V}_{i,k+1} \leq (\dim D(1, \theta))^{k-i} (\dim V(\theta))^{i+1} + (\dim D(1, \theta))^{k-i}
(\dim V(\theta))^{i}.
\end{gather*}
Together with~\eqref{dimd}, we see
\begin{gather*}
\dim \mathbf{V}_{i,k+1} \leq (\dim D(1,\theta))^{k+1-i} (\dim V(\theta))^{i}.
\end{gather*}
Now the proof of Theorem~\ref{MT} in this case follows by~\eqref{dimg}.
This completes the proof of Theorem~\ref{MT}.

Combining parts (1) and (2) of Theorem~\ref{MT}, we get Corollary~\ref{c1}.
Using part (2) of Theorem~\ref{MT} and Proposition~\ref{F}, we obtain Corollary~\ref{c2}.

\section{CV modules and truncated Weyl modules}\label{section5}

We start this section by recalling the def\/inition of CV modules given in~\cite{CV}.
For $\mathfrak{g}$ simply laced, we shall restate Theorem~\ref{MT} in terms of these modules.
At the end, we also discuss truncated Weyl modules.

{\bf 5.1.}~Given $\lambda\in P^{+}$, we say that $\pmb{\xi}=(\xi(\alpha))_{\alpha\in R^+}$ is a~$\lambda$-compatible
$|R^{+}|$-tuple of partitions, if
\begin{gather*}
\xi(\alpha)=\big(\xi(\alpha)_1\geq\dots \geq\xi(\alpha)_j\geq\dots\geq 0\big),
\qquad
|\xi(\alpha)|= \sum\limits_{j\geq 1}\xi(\alpha)_j = \langle \lambda, \alpha^{\vee} \rangle,
\qquad
\forall\, \alpha \in R^+.
\end{gather*}
\begin{definition}[see~\protect{\cite[\S~2]{CV}}] Let $\lambda\in P^{+}$ and $\pmb{\xi}$ be a~$\lambda$-compatible $|R^{+}|$-tuple of partitions.
The Chari--Venkatesh module or CV module $V(\pmb\xi)$ is the graded quotient of $W(\lambda)$ by the submodule generated
by the following set
\begin{gather*}
\bigg\{\mathbf{x}^{-}_{\alpha}(r,s)w_{\lambda}: \alpha\in R^+, s, r \in \mathbb{N}~\textup{such that}~s+r\geq
1+rk+\sum\limits_{j\geq k+1} \xi(\alpha)_j~\textup{for some}~k\in \mathbb{N}\bigg\}.
\end{gather*}
\end{definition}

The following lemma (implicit in the proof of Theorem~1 of~\cite{CV}) is useful in understanding CV modules.
\begin{lemma}
\label{c-v}
Let $\lambda\in P^{+}$, $r\in\mathbb{N}$ and $\pmb\xi=(\xi(\alpha))_{\alpha\in R^+}$ a~$\lambda$-compatible
$|R^{+}|$-tuple of partitions.
If $r\geq \xi(\alpha)_1$ then $\mathbf{x}^{-}_{\alpha}(r,s)w_{\lambda}=0$ in $W(\lambda)$, for all $\alpha\in R^+$,
$s, k \in \mathbb{N}$, $s+r\geq 1+rk+\sum\limits_{j\geq k+1} \xi(\alpha)_j$.
\end{lemma}
\begin{proof}
Let $\alpha\in R^+$ and $s, k \in \mathbb{N}$ such that $s+r\geq 1+rk+\sum\limits_{j\geq k+1} \xi(\alpha)_j$. Given
$r\geq \xi(\alpha)_1$, it follows that $s+r \geq 1 + \sum\limits_{j\geq 1} \xi(\alpha)_j=1+\langle \lambda,
\alpha^{\vee} \rangle$.
Now the proof follows by using Lemma~\ref{garland} and~\eqref{w2}.
\end{proof}

For $\lambda\in P^{+}$, we associate two~$\lambda$-compatible $|R^{+}|$-tuple of partitions as follows:
\begin{gather*}
\{\lambda\}:=\big(\big(\langle \lambda, \alpha^{\vee} \rangle\big)\big)_{\alpha\in R^{+}},
\qquad
\pmb\xi(\lambda):=\big(\big(1^{\langle \lambda, \alpha^{\vee} \rangle}\big)\big)_{\alpha\in R^+}.
\end{gather*}
Each partition of $\{\lambda\}$ has at most one part, and each part of each partition of $\pmb\xi(\lambda)$ is 1.
The CV modules corresponding to these two, have nice descriptions, which we record below for later use;
\begin{gather}
\label{e1}
V(\{\lambda\})\cong_{\mathfrak{g}[t]} \ev_0 V(\lambda),
\qquad
V(\pmb\xi(\lambda)) \cong_{\mathfrak{g}[t]} W(\lambda).
\end{gather}
The f\/irst isomorphism follows by taking $s=r=k=1$ in the def\/inition of the CV mo\-du\-le~$V(\{\lambda\})$ and the second
isomorphism follows from Lemma~\ref{c-v}.

{\bf 5.2.}~Given $k\geq 1$ and $0\leq i \leq k $, we def\/ine the following $|R^{+}|$-tuple of partitions:
\begin{alignat*}{3}
&\pmb\xi_{i}^{-}:=(\xi^{-}_{i}(\alpha))_{\alpha\in R^{+}},
\qquad&&
\text{where}
\quad
\xi^{-}_{i}(\alpha)=
\begin{cases}
\big(1^{\langle k\theta, \alpha^{\vee} \rangle}\big), & \alpha\neq \theta,
\\
\big(2^{i}, 1^{2(k-i)}\big), & \alpha=\theta,
\end{cases}&
\\
& \pmb\xi_{i}:=\big(\xi_{i}(\alpha)\big)_{\alpha\in R^{+}},
\qquad &&
\text{where}
\quad
\xi_{i}(\alpha)=
\begin{cases}
\big(1^{\langle (k+1)\theta, \alpha^{\vee} \rangle}\big), & \alpha\neq \theta,
\\
\big(2^{i}, 1^{2(k+1-i)}\big), & \alpha=\theta,
\end{cases}&
\\
& \pmb\xi_{i}^{+}:=\big(\xi^{+}_{i}(\alpha)\big)_{\alpha\in R^{+}},
\qquad &&
\text{where}
\quad
\xi^{+}_{i}(\alpha)=
\begin{cases}
\big(1^{\langle (k+1)\theta, \alpha^{\vee} \rangle}\big), & \alpha\neq \theta,
\\
\big(2^{i+1}, 1^{2(k-i)}\big), & \alpha=\theta.
\end{cases}&
\end{alignat*}

For $\mathfrak{g}$ simply laced, we can restate Theorem~\ref{MT} in terms of CV modules as follows:
\begin{theorem}
\label{T2}
Assume that $\mathfrak{g}$ is simply laced.
Given $k\geq 1$ and $0\leq i \leq k $, we have the following:
\begin{enumerate}\itemsep=0pt
\item[$1)$] a~short exact sequence of $\mathfrak{g}[t]$-modules,
\begin{gather*}
0 \rightarrow \tau_{2k+1-i}V(\pmb\xi^{-}_{i}) \rightarrow V(\pmb\xi_{i}) \rightarrow V(\pmb\xi^{+}_{i}) \rightarrow 0;
\end{gather*}
\item[$2)$] an isomorphism of $\mathfrak{g}[t]$-modules,
\begin{gather*}
V(\pmb\xi_{i}) \simeq V(\pmb\xi(\theta))^{* (k+1-i)} * V(\{\theta\})^{*i}.
\end{gather*}
\end{enumerate}
\end{theorem}
\begin{proof}
This follows from Theorem~\ref{MT}, by using Lemma~\ref{c-v} and~\eqref{e1}.
\end{proof}

{\bf 5.3.}~For $n\geq 1$, we def\/ine $\mathcal{A}_n= \mathbb{C}[t]/(t^{n})$.
The {\em truncated current algebra} $\mathfrak{g}\otimes \mathcal{A}_n$, can be thought of as the graded quotient of the
current algebra $\mathfrak{g}[t]$:
\begin{gather*}
\mathfrak{g}\otimes \mathcal{A}_n\cong\mathfrak{g}[t]/\big(\mathfrak{g}\otimes t^n \mathbb{C}[t]\big).
\end{gather*}
Let $k\geq1$.
The local Weyl module $W_{\mathcal{A}_n}(k\theta)$ for the truncated current algebra $\mathfrak{g}\otimes \mathcal{A}_n$
is def\/ined in~\cite{CFK}, and we call it the {\em truncated Weyl module}.
It is easy to see that $W_{\mathcal{A}_n}(k\theta)$ naturally becomes a~$\mathfrak{g}[t]$-module and the following is an
isomorphism of $\mathfrak{g}[t]$-modules,
\begin{gather}
\label{trul}
W_{\mathcal{A}_n}(k\theta) \cong  W(k\theta)/\langle (x^{-}_{\theta}\otimes t^{n})  w_{k\theta}\rangle.
\end{gather}
Now Corollary~\ref{truncated} is immediate from Corollary~\ref{c2}, by using~\eqref{w2'} and~\eqref{trul}.

\subsection*{Acknowledgements}

The author thanks Vyjayanthi Chari, K.N.~Raghavan and S.~Viswanath for many helpful discussions and encouragement.
Part of this work was done when the author was visiting the Centre de Recherche Mathematique (CRM), Montreal, Canada,
during the thematic semester on New Directions in Lie Theory.
The author acknowledges the hospitality and f\/inancial support extended to him by CRM.
The author also thanks the anonymous referees for their valuable comments, due to which the paper is much improved.
The author acknowledges support from CSIR under the SPM Fellowship scheme.

\pdfbookmark[1]{References}{ref}
\LastPageEnding


\begin{thebibliography}{99}
\footnotesize \itemsep=0pt

\bibitem{CFK}
Chari V., Fourier G., Khandai T., A categorical approach to {W}eyl modules,
  \href{http://dx.doi.org/10.1007/s00031-010-9090-9}{\textit{Transform. Groups}} \textbf{15} (2010), 517--549, \href{http://arxiv.org/abs/0906.2014}{arXiv:0906.2014}.

\bibitem{CP}
Chari V., Pressley A., Weyl modules for classical and quantum af\/fine algebras,
  \href{http://dx.doi.org/10.1090/S1088-4165-01-00115-7}{\textit{Represent. Theory}} \textbf{5} (2001), 191--223,
  \href{http://arxiv.org/abs/math.QA/0004174}{math.QA/0004174}.

\bibitem{CSVW}
Chari V., Shereen P., Venkatesh R., Wand J., A Steinberg type decomposition
  theorem for higher level Demazure modules, \href{http://arxiv.org/abs/1408.4090}{arXiv:1408.4090}.

\bibitem{CV}
Chari V., Venkatesh R., Demazure modules, fusion products, and $Q$-systems,
  \href{http://dx.doi.org/10.1007/s00220-014-2175-x}{\textit{Comm. Math. Phys.}}, {t}o appear, \href{http://arxiv.org/abs/1305.2523}{arXiv:1305.2523}.

\bibitem{FL}
Feigin B., Loktev S., On generalized {K}ostka polynomials and the quantum
  {V}erlinde rule, in Dif\/ferential Topology, Inf\/inite-Dimensional {L}ie
  Algebras, and Applications, \textit{Amer. Math. Soc. Transl. Ser.~2}, Vol.~194, Amer. Math. Soc., Providence, RI, 1999, 61--79,
  \href{http://arxiv.org/abs/math.QA/9812093}{math.QA/9812093}.

\bibitem{F}
Feigin E., The {PBW} f\/iltration, {D}emazure modules and toroidal current
  algebras, \href{http://dx.doi.org/10.3842/SIGMA.2008.070}{\textit{SIGMA}} \textbf{4} (2008), 070, 21~pages,
  \href{http://arxiv.org/abs/0806.4851}{arXiv:0806.4851}.

\bibitem{FoL}
Fourier G., Littelmann P., Weyl modules, {D}emazure modules, {KR}-modules,
  crystals, fusion products and limit constructions, \href{http://dx.doi.org/10.1016/j.aim.2006.09.002}{\textit{Adv. Math.}}
  \textbf{211} (2007), 566--593, \href{http://arxiv.org/abs/math.RT/0509276}{math.RT/0509276}.

\bibitem{G}
Garland H., The arithmetic theory of loop algebras, \href{http://dx.doi.org/10.1016/0021-8693(78)90294-6}{\textit{J.~Algebra}}
  \textbf{53} (1978), 480--551.

\bibitem{Naoi}
Naoi K., Weyl modules, {D}emazure modules and f\/inite crystals for non-simply
  laced type, \href{http://dx.doi.org/10.1016/j.aim.2011.10.005}{\textit{Adv. Math.}} \textbf{229} (2012), 875--934,
  \href{http://arxiv.org/abs/1012.5480}{arXiv:1012.5480}.

\bibitem{V}
Venkatesh R., Fusion product structure of Demazure modules, \href{http://dx.doi.org/10.1007/s10468-014-9495-6}{\textit{Algebr.
  Represent. Theory}}, {t}o appear, \href{http://arxiv.org/abs/1311.2224}{arXiv:1311.2224}.

\end{thebibliography}
\end{document}